\documentclass[11pt]{article}

\def\mathbb{\Bbb}

\usepackage{amsfonts,latexsym,amstext}
\usepackage{amsmath}
\usepackage{amssymb}

\usepackage{comment}

\newtheorem{theorem}{Theorem}[section]
\newtheorem{lem}[theorem]{Lemma}
\newtheorem{proposition}[theorem]{Proposition}
\newtheorem{definition}{Definition}[section]

\newtheorem{hypothesis}[theorem]{Hypothesis}

\newtheorem{remark}[theorem]{Remark}

\makeatletter
\@addtoreset{equation}{section}

\makeatother

\def\qed{{\hfill\hbox{\enspace${ \square}$}} \smallskip}
\def\sqr#1#2{{\vcenter{\vbox{\hrule height .#2pt \hbox{\vrule
 width .#2pt height#1pt \kern#1pt \vrule
width .#2pt} \hrule height .#2pt}}}}
\def\square{\mathchoice\sqr54\sqr54\sqr{4.1}3\sqr{3.5}3}

\def\ds{\begin{displaystyle}}
\def\eds{\end{displaystyle}}
\def\dis{\displaystyle }
\def\<{\langle }
\def\>{\rangle }

\newcommand{\vsc}{\vskip 5mm}

\def\R{\mathbb R}

\def\E{\mathbb E}
\def\P{\mathbb P}

\def\L{\mathbb L}

\def\cala{{\cal A}}

\def\calf{{\cal F}}

\def\calk{{\cal K}}

\def\calp{{\cal P}}

\def\calu{{\cal U}}

\newcommand{\sign}{\mbox{\rm sign}}

\newcommand{\La}{\Lambda}

\newcommand{\si}{\sigma}

\newcommand{\be}{\beta}
\newcommand{\ep}{\varepsilon}

\newcommand{\De}{\Delta}
\newcommand{\om}{\omega}
\newcommand{\Om}{\Omega}
\newcommand{\ze}{\zeta}

\newcommand{\aaa}{\mathcal{A}}
\newcommand{\ba}{\mathcal{B}}

\newcommand{\f}{\mathcal{F}}
\newcommand{\g}{\mathcal{G}}
\newcommand{\h}{\mathcal{H}}

\newcommand{\p}{\mathcal{P}}

\newcommand{\s}{\mathcal{S}}

\newcommand{\ua}{\mathcal{U}}

\newcommand{\bee}{\begin{equation}}
\newcommand{\eee}{\end{equation}}
\newcommand{\bea}{\begin{eqnarray}}
\newcommand{\eea}{\end{eqnarray}}
\newcommand{\bean}{\begin{eqnarray*}}
\newcommand{\eean}{\end{eqnarray*}}

\title{$L^p$ solution of backward stochastic differential equations driven by a marked point process}

\date{}
\author{Fulvia Confortola
\\Politecnico di Milano,
Dipartimento di Matematica\\
piazza Leonardo da Vinci 32, 20133 Milano, Italy\\
e-mail: fulvia.confortola@polimi.it}

\begin{document}

\maketitle

\begin{abstract}\noindent
\\
We  obtain existence and uniqueness in  $L^p$, $p>1$ of the solutions of a backward stochastic differential equations (BSDEs
for short)
driven by a marked point process, on a bounded interval.
   We show that the  solution of the BSDE can be approximated by a finite system  of deterministic differential equations.
As application we address an
optimal control problems for point processes of general non-Markovian type and show that BSDEs
can be used to prove existence of an optimal control and to represent the value function.

\end{abstract}

\section{Introduction}\label{s-Intro}
\setcounter{equation}{0}
In this paper we consider a backward stochastic differential equation
(BSDE for short in the remaining) driven by a random measure, without diffusion
part, on a finite time interval, of the following form:
\bee\label{Intro1}
Y_t+\int_t^T\int_K Z(s,x)\,\mu(ds,dx)
=\xi+\int_t^T\int_K f(s,x,Y_{s-},Z_s(\cdot))\,\nu(ds,dx).
\eee
 $\mu$ is the counting measure corresponding to a non-explosive marked point process $(S_n,X_n)_{n\ge 1}$, where
$(S_n)$ is an increasing sequence of random times and $(X_n)$ a sequence of random
variables in the state (or  mark) space $K$.  $\nu$ is the predictable compensator of the measure $\mu$. The generator $f$ and the final condition $\xi$ are given. The unknown process is a pair
 $(Y_t,Z_t(\cdot))$, where $Y$ is a real adapted c\`adl\`ag process and
$\{Z_t(x),\, t\in [0,T],\,x\in K\}$ is a predictable random field.

The BSDEs have been introduced by Pardoux and Peng \cite{PP90}. In this first paper, and in most of the subsequent ones, the driving
term is a Brownian motion. Since then, there has been an increasing interest
for this subject: these equations have
a wide range of applications in  to various fields of
stochastic analysis, including probabilistic techniques in partial
differential equations, stochastic optimal
control, mathematical finance (see e.g. El Karoui, Peng and
Quenez \cite{EPQ}). Recently, BSDEs driven by a Brownian motion and a random measure
have also been considered due their utility in the study of stochastic maximum principle, partial differential equations
of nonlocal type, quasi-variational inequalities, impulse control and stochastic problem in stochastic finance, see e.g. Buckdahn and Pardoux
\cite{BP94}, Tang and Li \cite{tali}, Barles, Buckdahn and
Pardoux \cite{BBP}, Xia \cite{xia}, Becherer \cite{Bech}, Cr\'epey and
Matoussi \cite{CM-2008}, or Carbone, Ferrario and Santacroce
\cite{CaFeSa} among many others.

In spite of the large literature devoted to BSDEs  with driving term  continuous or
continuous-plus-jumps there are relatively few results on the case of a driving term which is purely discontinuous.
We cite Shen and Elliott \cite{SE-2011}
for the particularly simple ``one-jump'' case, or Cohen and Elliott
\cite{CoEl}-\cite{CoEl-a} and Cohen and Szpruch \cite{CoSz} for BSDEs
associated to Markov chains.

In \cite{CoFu-m} and \cite{BaCo} are considered BSDEs driven by more general random measures, related to a Markov and semi-Markov process respectively, in connection with optimal control problem. The more general non-Markovian case is studied in \cite{CoFu-mpp}. Here the authors, relying upon the martingale representation theorem, prove existence, uniqueness and continuous dependence on the data of the solution of equation \eqref{Intro1} in suitable weighted $L^2$ spaces. They require  a $L^2$ summability condition on the data $\xi$ and $f$ and a Lipschitz condition on the generator $f$.

The first aim of this paper is  to extend the results contained in \cite{CoFu-mpp}  and to develop a $L^p$-theory for $p>1$ for this class of BSEDs.
 The basic hypothesis on the generator $f$ is an uniform Lipschitz
condition (see Hypothesis \ref{hyp:BSDE} for precise
statements). In order to solve the equation, beside measurability assumptions,
  we require, for $p>1$
 the $L^p$ summability condition
 $$
\E \, [e^{\beta A_T}|\xi|^p] + \E \int_0^T \int_K e^{\beta A_t}|f(\om,t,0,0)|^p \nu(dt,dx)
  < \infty,
 $$
 to hold for  a suitable  $\beta$. To prove existence and uniqueness for the BSDE we require -as in \cite{CoFu-mpp}- only that the jump times $S_n$ are totally inaccessible (see \textbf{Assumption }$\bf{A_1}$ below).

The results are stated in the
case of a scalar equation, but the extension to the vector-valued case is immediate. They are presented in Section 3, after an introductory
section devoted to notation and preliminaries.

We recall that the $L^1$ theory for the solutions of the equation \eqref{Intro1} is developed in \cite{CFJ}. Our $L^p$ assumptions are not in general comparable with the $L^1$ assumptions which involve suitable doubly weighted spaces. Moreover we require only that the point process is \emph{quasi-left continuous} (\textbf{Assumption }$\bf{A_1}$) and, using the martingale representation theorem and fixed point arguments, we do not need the technical \textbf{Assumption}  $\bf{A_2}$ (see Subsection \ref{subsect: ode}).

The second purpose of the paper is presented
in  Section \ref{S4}, where  we illustrate an approximation scheme to solve the equation \eqref{Intro1}.

\noindent  We show that the solution to \eqref{Intro1} can be obtained as limit of a sequence of approximating BSDEs driven by a random measure with a finite number of jumps.
We use  the \emph{a priori} estimates to obtain an error estimate for this
approximation (see Proposition \ref{approxy,z}). It involves  the distribution of the last jump time and in particular cases it can be easily computed.

If we add now  \textbf{Assumption} $\bf{A_2}$, we can  replace the approximating BSDE with a system of finite deterministic ordinary differential equation, using the result in \cite{CFJ}. The method to reduce the BSDE to a sequence of ODEs  has been  used also for a BSDE driven by a Brownian
motion plus a Poisson process, see e.g. Kharroubi and Lim [17].
In recent years there has been much interest in numerical approximation of the solution to the BSDEs, in the context of diffusion processes.
Our results might be used for similar methods in the framework of pure jump
processes as well.

In Section \ref{S5} we address an optimal control problem for a marked process, formulated in a classical way, with the BSDEs approach. We extend  the results on optimal control problem in  \cite{ CoFu-mpp}, assuming, for $p>1$, $L^p$ summability conditions on the data of the problem.

\section{The setting}
Let $T \in (0, \infty)$ be a fixed time horizon.
Let $(\Omega,\calf, \P)$ be a complete probability space
and $(K,\calk)$ a Lusin space.
Assume we have a  non-explosive multivariate point process
(also called marked point process) on $[0,T]\times K$: this is a sequence $(S_n,X_n)_{n\ge 1}$ of
random variables  with
distinct times of occurrence $S_n$ and with marks $X_n$. $S_n$ taking values in $(0,T]\cup\{\infty\}$
and $X_n$ in $K$. We set $S_0=0$
and we assume, $\P$-a.s., $S_1>0$; if $S_n \leq T$ then $ S_n<S_{n+1}$ and $S_n\leq
S_{n+1}$ everywhere; $\Om=\cup\{S_n>T\}$. Note that the ``mark'' $X_n$ is
relevant on the
set $\{S_n\leq T\}$ only,  but it is convenient to have it defined on the
whole set $\Om$, and without restriction we may assume that $X_n=\De$ when
$S_n=\infty$, where $\De$ is a distinguished point in $K$.

The multivariate point process can be
viewed as a random measure of the form
\bee\label{S1}
\mu(dt,dx)=\sum_{n\geq1:\,S_n\leq T}\ep_{(S_n,X_n)}(dt,dx).
\eee
where $\ep_{(t,x)}$ denotes the Dirac measure.

We denote by $(\f_t)_{t\geq0}$ the filtration generated by the point process,
which is the smallest filtration for which each $S_n$ is a stopping time and
$X_n$ is $\f_{S_n}$-measurable. As we will see, the special structure of this
filtration plays a fundamental role in all what follows. We let
$\p$ be the predictable $\si$-field on $\Om\times[0,T]$, and for any auxiliary
measurable space $(G,\g)$ a function on the product $\Om\times[0,T]\times G$
which is measurable with respect to $\p\otimes\g$ is called
{\em predictable}.

We denote by $\nu$ the predictable compensator of the measure $\mu$,
relative to the filtration $(\f_t)$. The
measure $\nu$ admits the disintegration:
\bee\label{S2}
\nu(\om,dt,dx)~=~dA_t(\om)\,\phi_{\om,t}(dx),
\eee
where $A$ is an increasing c\`adl\`ag predictable process
starting at $A_0=0$, which is also the predictable compensator of the
univariate point process
\bee\label{S3}
N_t~=~\mu([0,t]\times E)~=~\sum_{n\geq1}1_{\{S_n\leq t\}};
\eee
$\phi$ is a transition probability from $(\Om\times[0,T],\p)$ into
$(K,\mathcal{K})$ and verifies the following equality
\bee \label{phi}
\E \int_{0}^\infty \int_K H(t,x)\; \mu(dt\, ,dx)=\E \int_{0}^\infty \int_K H(t,x)\;\phi_t(dx)\, dA_t
\eee
for every nonnegative $H_t(\omega,y)$, $\calp\otimes \calk$-measurable.

The following assumption
will hold throughout:
 $ $

\noindent \textbf{Assumption $\bf{A_1}$}: \rm The process $A$ is continuous (equivalently: the
jump times $S_n$ are totally inaccessible).

This condition amounts to the quasi-left continuity of $N$.

\section{The backward equation} \label{sec-backward}

We denote by $\ba(K)$ the set of all Borel
functions on $K$; if $Z$ is a measurable function on $\Om\times[0,T]\times K$,
we write $Z_{\om,t}(x)=Z(\om,t,x)$, so each $Z_{\om,t}$, often abbreviated
as $Z_t$ or $Z_t(\cdot)$, is an element of $\ba(K)$.

In the following we will consider the backward stochastic differential equation
\bee\label{S11}
Y_t+\int_t^T\int_K Z(s,x)\,\mu(ds,dx)
=\xi+\int_t^T\int_K f(s,x,Y_{s-},Z_s(\cdot))\,\nu(ds,dx),
\eee
where the generator $f$ and the final condition $\xi$ are given.

\begin{definition} \label{def:solBSDE} A {\em solution} is a pair $(Y,Z)$ consisting in an adapted c\`adl\`ag
process $Y$ and a predictable function $Z$ on $\Om\times[0,T]\times E$
satisfying $\int_0^T\int_E|Z(t,x)|\,\nu(ds,dx)<\infty$ a.s., such
that \eqref{S11} holds for all $t\in[0,T]$, outside a $\P$-null set.
\end{definition}
An other notion of \emph{solution} can be introduced by observing that
\eqref{S11} ca be rewritten as follows:
\bee\label{S15}
Y_t+\sum_{n\geq1}Z(S_n,X_n)\,1_{\{t<S_n\leq T\}}
=\xi+\int_{(t,T]}\int_Ef(s,x,Y_{s-},Z_s(\cdot))\,\nu(ds,dx).
\eee
Since $A$ is continuous, \eqref{S15} yields, outside a $\P$-null set:
\bee\label{S16}
\De Y_{S_n}=Z(S_n,X_n)~~\text{if $S_n\leq T$ and $n\geq1$},\qquad
Y~\text{is continuous outside $\{S_1,\cdots,S_n,\cdots\}$}.
\eee
In other words, $Y$ completely determines the predictable function
$Z$ outside a null set with respect to the measure
 $\P(d\om)\mu(\om,dt,dx)$, hence also
outside a $\P(d\om)\nu(\om,dt,dx)$-null set. Equivalently, if $(Y,Z)$ is a
solution and $Z'$ is another predictable function, then $(Y,Z')$ being another
solution is the same as having $Z'=Z$ outside a $\P(d\om)\mu(\om,dt,
dx)$-null set, and the same as having $Z'=Z$ outside a $\P(d\om)
\nu(\om,dt,dx)$-null set.

Hence it is possible to define
a {\em solution} to \eqref{S11} an adapted c\`adl\`ag process $Y$ for which there exists
a predictable function $Z$ satisfying $$\int_0^T\int_E|Z(s,x)|\,\nu(ds,dx)
<\infty \mbox{ a.s.},$$ such that the pair $(Y,Z)$ satisfies \eqref{S11} for all
$t\in[0,T]$, outside a $\P$-null set. Then, {\em uniqueness} of the solution
means that, for any two solutions $Y$ and $Y'$ we have $Y_t=Y'_t$ for all
$t\in[0,T]$, outside a $\P$-null set.

\subsection{The $L^p$ theory.} We introduce the Banach space $L^p_{\beta}$, depending on a parameter $\beta>0$,
of equivalence classes of pairs of processes $(Y,Z)$ on $[0,T]$
such that $Y$ is progressive, $Z$ is predictable and the norm
$$
\|(Y,Z)\|=\E\left[
    \int_0^T\int_K(|Y_s |^p+|Z(t,x) |^p)\,e^{\beta A_s}\, \nu(dt,dx)  \right]
$$
is finite. Elements of $L^p_{\beta}$ are identified up to almost sure
 equality with respect to the measure
$\P(d\omega)\,d\nu(dt,dx)$, i.e. when the norm of their difference is zero.
We sometimes identify processes $(Y,Z)$ with their equivalence classes
in the usual way.

Let us consider the following assumptions on the data $\xi$ and  $f$:
\begin{hypothesis}\label{hyp:BSDE}
\begin{enumerate}
\item The final condition $\xi:\Omega\to\R$ is $\calf_T$-measurable and
$\E \, e^{\beta A_T}|\xi|^p < \infty$.

\item $f$ is a real-valued function on $\Om\times[0,T]\times K\times\R
\times\ba(K)$, such that
\begin{itemize}
\item[(i)] for any predictable
function $Z$ on $\Om\times[0,T]\times K$ the mapping
\begin{equation}\label{fmisurabile}
(\om,t,x,y)  \mapsto f(\om,t,x,y,Z_{\om,t}(\cdot))
\end{equation} is predictable;
\item[(ii)] there exist $L \geq0$, $L' \geq 0$ such that for every $\omega \in \Om$, $t \in [0,T]$, $x \in K$, $y,y' \in \mathbb{R}$, $\ze', \ze \in \ba(K)$  we have
    \bee \label{hyp:f-lip}
    \begin{array}{l}
    |f(\om,t,x,y',\ze)-f(\om,t,x,y,\ze)|\leq L'|y'-y|\\
    \begin{array}{lll} \int_K|f(\om,t,x,y,\ze)-f(\om,t,x,y,\ze')|\phi_{\om,t}(dx) & \leq & L\left(\int_K|\ze'(v)-\ze(v)|^p\,\phi_{\om,t}(dv)\right)^{1/p}\\
                                                                                  & = & L|\ze- \ze'|_{L^p(\phi)}
\end{array}
\end{array}
\eee
where $\phi_{\om,t}$ are the measures occurring in
\eqref{S2};
\item[(iii)] We have
\begin{equation} \label{Lp-f0}
\E\int_0^T\int_K e^{\beta A_t}|f(t,x,0,0)|^p\,\nu(dt,dx)<\infty \quad a.s.
\end{equation}
\end{itemize}
\end{enumerate}
\end{hypothesis}

The measurability condition imposed on the generator is slightly involved, but it is verified when we deal with a BSDE in order to solve an optimal control problem driven by a multivariate point processes. In this framework  the suitable formulation is the following one
\bee\label{S9}
Y_t+\int_t^T\int_KZ(s,x)\,\mu(ds,dx)
=\xi+\int_t^T\bar{f}(.,s,Y_{s-},\eta_sZ_s)\,dA_s,
\eee
where
\bee\label{S10}
\begin{array}{c}
\text{$\eta_{\om,t}$ is a real-valued map on $\ba(K)$, with}~~
|\eta_{\om,t}\ze-\eta_{\om,t}\ze'|\leq \int_K|\ze'(v)-\ze(v)|
\,\phi_{\om,t}(dv),\\
\text{$Z$ predictable on}~ \Om\times[0,T]\times K~~\Rightarrow~~
\text{the process $(\om,t)\mapsto \eta_{\om,t}Z_{\om,t}$ is predictable,}
\end{array}
\eee
$\bar{f}$ is a predictable function on
$\Om\times[0,T]\times\R\times\R$ satisfying:
\bee
\begin{array}{c}\label{fbar}
|\bar{f}(t,y',z')-\bar{f}(t,y,z)|\leq L'|y'-y|+L|z'-z|\\
\E \int_0^T \int_K|\bar{f}(t,x,0,0)|^p e^{\beta A_t}\,\nu(dt,dx)<\infty
\end{array}
\eee

The equation \eqref{S9} reduces to \eqref{S11} upon taking
\bee\label{S14}
f(\om,s,x,y,\ze)=\bar{f}(\om,s,y,\eta_{\om,s}\ze),
\eee
and \eqref{fbar} for $\bar{f}$ plus \eqref{S10} for $\eta_{\om,t}$
yield, with the H\"{o}lder inequality, the \eqref{hyp:f-lip} for $f$.

\subsubsection{A priori estimate}

In this section, we provide some {\em a priori} estimates for the solutions
of Equation \eqref{S11}.

We start with a Lemma of Ito type.

\begin{lem} \label{lem:ito} Let $\beta\in\R$. If $(Y,Z) \in L^p_{\beta}$
is a solution of \eqref{S11} we have almost surely

\begin{equation}\label{P1}
\begin{array}{c}
|Y_t|^pe^{\beta A_t}+\int_t^T\int_K\big(|Y_{s-}+Z(s,y)|^p
-|Y_{s-}|^p\big)\,e^{\beta A_s}\,\mu(ds,dy)+\beta\int_t^T|Y_s|^p
e^{\beta A_s}\,dA_s\\[.3cm]
=|\xi|^p\,e^{\beta A_T}+\int_t^T\int_K p|Y_{s-}|^{p-1}\sign(Y_s-)
\,f(s,y,Y_s,Z_s(\cdot))\,e^{\beta A_s}\,\nu(ds,dy).
\end{array}
\end{equation}
\end{lem}

\noindent \textbf{Proof}.  Letting $U_t$ and $V_t$ be the left and right sides of
%\eqref{P1},
and since these processes are c\`adl\`ag, and continuous
outside the $S_n$'s, and $U_T=V_T$, it suffices to check that outside a null
set we have $\Delta U_{S_n}=
\Delta V_{S_n}$ and also $U_t-U_s=V_t-V_s$ if $S_n\leq t<s<S_{n+1}\wedge T$, for
all $n\geq0$. The first property is obvious because $\Delta Y_{S_n}=Z(S_n,X_n)$
a.s. and $A$ is continuous.
The second property follows from
$Y_t-Y_s=\int_t^s\int_K f(v,y,Y_v,Z_v(\cdot))\,\nu(dv,dy)$,
implying $|Y_t|^p-|Y_s|^p=\int_t^s\int_K p|Y_r|^{p-1}\sign(Y_s-)\,f(r,y,Y_r,
Z_r(\cdot))\,\nu(dr,dx)$
 plus a standard change
of variables formula.\qed

\begin{comment}

\begin{lem} Let $\beta\in\R$. If $(Y,Z)$
is a solution of  we have almost surely

\begin{equation}\label{P1}
\begin{array}{c}
|Y_t|^pe^{\beta A_t}+\int_t^T\int_K\big(|Y_{s-}+Z(s,y)|^p
-|Y_{s-}|^p\big)\,e^{\beta A_s}\,\mu(ds,dy)+\beta\int_t^T|Y_s|^p
e^{\beta A_s}\,dA_s\\[.3cm]
=|\xi|^p\,e^{\beta A_T}+\int_t^T p|Y_{s-}|^{p-1}
\,f_s\,e^{\beta A_s}\,\nu(ds,dy).
\end{array}
\end{equation}
\end{lem}

\noindent \textbf{Proof}.  Letting $U_t$ and $V_t$ be the left and right sides of
%\eqref{P1},
and since these processes are c\`adl\`ag, and continuous
outside the $S_n$'s, and $U_T=V_T$, it suffices to check that outside a null
set we have $\Delta U_{S_n}=
\Delta V_{S_n}$ and also $U_t-U_s=V_t-V_s$ if $S_n\leq t<s<S_{n+1}\wedge T$, for
all $n\geq0$. The first property is obvious because $\Delta Y_{S_n}=Z(S_n,X_n)$
a.s. and $A$ is continuous.
The second property follows from
$Y_t-Y_s=\int_t^s\int_K f_s\,\nu(dv,dy)$,
implying $|Y_t|^p-|Y_s|^p=\int_t^s\int_K p|Y_r|^{p-1}\,f_s\,\nu(dr,dx)$
 plus a standard change
of variables formula.\qed
\end{comment}

$ $

Next we prove the following useful a priori estimates:
\begin{lem}\label{LP2} For every $\epsilon >0$ let  $C_{\epsilon}=\left(1-\left(\frac{1}{1+ \epsilon}\right)^{\frac{1}{p-1}}\right)^{1-p}$.

\noindent Suppose that Hypothesis \ref{hyp:BSDE} holds with  $\beta> 1 + \frac{c_{\epsilon}}{1+ \epsilon}+ pL'+(p-1)(L^p(1+ \epsilon))^{\frac{1}{p-1}}$. Then there exist two constants
$C_1$ and $C_2$ only depending on $(\beta,p,L,L',\epsilon)$, such
that any pair $(Y,Z)$ in $L^p_{\beta}$
which solves \eqref{S11}
satisfies
\bee\label{P3}
\E|Y_ t|^p e^{\be A_t}\leq C_1\E\Big(|\xi|^pe^{\be A_T}+
\int_t^T\int_K |f(s,x,0,0)|^p\,e^{\be A_s}\,\nu(ds,dy)\Big)
\eee
\bee\label{P4}
\E\int_0^T\int_K(|Y_{s}|^p+ |Z(s,x)|^p)\,e^{\be A_s}\,\nu(ds,dx)\leq
C_2\,\E\Big(|\xi|^p e^{\be A_T}+\int_0^T\int_K
|f(s,x,0,0)|^p\,e^{\be A_s}\,\nu(ds,dy)\Big).
\eee
\end{lem}

\noindent \textbf{Proof. } We have
$|Y_{s-}+Z(s,x)|^p-|Y_{s-}|^p\geq \frac{1}{1+ \epsilon}|Z(s,x)|^p- (1 + \frac{c_{\epsilon}}{1+ \epsilon})|Y_{s-}|^p$, hence
\eqref{P1}, and the fact that
$\phi_{t,\om}(K)=1$ yield
almost surely,
\bee\label{P5}
\begin{array}{lll}
\dis |Y_t|^pe^{\be A_t}&+&\frac{1}{1+ \epsilon}\int_t^T\int_K|Z(s,x)|^p\,e^{\be A_s}
\,\mu(ds,dx)
+\be \int_t^T|Y_s|^pe^{\be A_s}\,dA_s \\
\dis &\leq&|\xi|^pe^{\be A_T}+
(1 + \frac{c_{\epsilon}}{1+ \epsilon})\int_t^T|Y_{s-}|^p\,e^{\be A_s} \,dN_s\\
\dis &+& p\int_t^T\int_K |Y_s|^{p-1}\,|f(s,x,Y_s,Z_s(\cdot)|\,
e^{\be A_s}\nu(ds,dx).
\end{array}
\eee

Taking the expectation in \eqref{P5} yields
\bee
\begin{array}{l}\label{stima1}
 \dis \E|Y_t|^p\,e^{\be A_t}+ \E\int_t^T\int_K\Big(\frac{1}{1+ \epsilon}|Z(s,x)|^p
+\left(\be - 1 - \frac{c_{\epsilon}}{1+ \epsilon} \right)|Y_s|^p\Big)\,e^{\be A_s}\,\nu(ds,dx)\\
\dis \quad\leq \E\Big(|\xi|^pe^{\be A_T}+
\int_t^T\!\!\int_K p|Y_s|^{p-1}\,|f(s,x,Y_s,Z_s(\cdot)|\,
e^{\be A_s}\,\nu(ds,dx)\Big)
\end{array}
\eee
From the Lipschitz condition of $f$ and elementary inequalities it follows that
$$
\begin{array}{l}
\dis \E|Y_t|^p\,e^{\be A_t}+ \E\int_t^T\int_K\Big(\frac{1}{1+ \epsilon}|Z(s,x)|^p
+\left(\be - 1 - \frac{c_{\epsilon}}{1+ \epsilon} \right)|Y_s|^p\Big)\,e^{\be A_s}\,\nu(ds,dx)\\
\dis \quad\leq \E\Big(|\xi|^pe^{\be A_T}+
\int_t^T\!\!\int_K p|Y_s|^{p-1}\,[|f(s,x,0,0)|+L'|Y_s |] \, e^{\be A_s}\,\nu(ds,dx) \Big)\\
 \dis \quad + \E\Big(\int_t^T p|Y_s|^{p-1}\, L|Z_s|_{L^p(\phi)} e^{\be A_s}\,dA_s\\
\quad\leq \E|\xi|^pe^{\be A_T}+\E\int_t^T\!\!\int_K \Big(\frac{p-1}{\gamma}\Big)^{p-1}|f(s,x,0,0)|^p \,e^{\be A_s}\,\nu(ds,dx)\\
\quad+\E\int_t^T\!\!\int_K\left[\left(pL' +(p-1)(L^p(\frac{1+ \epsilon}{\alpha}))^{\frac{1}{p-1}}+ \gamma \right) |Y_s|^p+ \frac{\alpha}{1+ \epsilon}|Z_s|^p \right]e^{\be A_s}\,\nu(ds,dx),
\end{array}$$
with $\alpha, \gamma >0$. If we choose $\alpha \in (0,1)$ and  $\gamma= \frac{1}{2}[\be- 1 - \frac{c_{\epsilon}}{1+ \epsilon}- pL'-(p-1)(L^p(\frac{1+ \epsilon}{\alpha}))^{\frac{1}{p-1}}]$, when $(Y,Z) \in L^p_{\beta}$, this implies almost surely,
$$\begin{array}{l}
\E|Y_t|^p\,e^{\be A_t}+ \frac{1- \alpha}{1+ \epsilon}\E\int_t^T\int_K|Z(s,x)|^p\,e^{\be A_s}\,\nu(ds,dx)\\
\quad + \E\int_t^T\int_K \frac{1}{2}\left(\be- 1 - \frac{c_{\epsilon}}{1+ \epsilon}- pL'-(p-1)(2L^p(1+ \epsilon))^{\frac{1}{p-1}}\right)|Y_s|^p\,\nu(ds,dx) \\ 
\leq \E\Big(|\xi|^pe^{\be A_T}+
\Big(\frac{p-1}{\gamma}\Big)^{p-1}\int_t^T\int_K|f(s,x,0,0)|^p\,e^{\be A_s}\,\nu(ds,dx)\Big),
\end{array}$$
giving us both \eqref{P3} and \eqref{P4}.\qed

%%%%%%%%%%%%%%%%%%%%%%%%%%%%%%%%%%%
\subsubsection{Existence and uniqueness}

In this section we will give an existence and uniqueness result for the  equation (\ref{S11}).
Our intent is to use the following integral representation theorem of marked
point process martingales (see e.g.  \cite{Da-art},\cite{Da-bo}).

\begin{theorem}\label{rappresentazione} \cite[Theorem 2.2]{CoFu-m}
Given $(t,x)\in [0,T]\times K$,
let $M$ be an $\mathcal{F}_{t}$-martingale on $[t,T]$ with respect
to $\P$. Then there exists a  predictable process $H $ on $\Omega \times [0,T]\times K$ with $\E \int_0^T \int_K |H_s(x)|\, \nu(ds , dx) <\infty$ such that
\begin{equation}\label{rapprmart}
    M_t= M_0 + \int_0^t \int_K H_s(x)\,( \mu(ds , dx)- \nu(ds,dx)),\qquad   t \in [0, T].
\end{equation}
\end{theorem}

\begin{theorem}\label{Thm:existence!bsde}
Suppose that  $ \xi$  and the generator $f$  satisfy Hypothesis \ref{hyp:BSDE}  and $\beta>1 + \frac{c_{\epsilon}}{1+ \epsilon}+ pL'+(p-1)((L+1)^p(1+ \epsilon))^{\frac{1}{p-1}}$. Then there exists a unique pair $(Y,Z)$ in $L^p_{\beta}$ which solves the BSDE \eqref{S11}.
\end{theorem}

\noindent \textbf{Proof. }
We start rewriting equation \eqref{S11} in the following equivalent way
\bee\label{S11-eq}
Y_t+\int_t^T\int_K Z(s,x)\,(\mu(ds,dx) - \nu(ds,dx))
=\xi+\int_t^T\int_K f(s,x,Y_{s-},Z(s,x)) -Z(s,x) \,\nu(ds,dx)
\eee
\noindent This formulation singles
out the ``martingale increment" $\int_t^T\int_K Z(s,x)\,(\mu(ds,dx) - \nu(ds,dx))$ and allows us to use Theorem \ref{rappresentazione}.
We define a mapping $\Phi$ from $L^p_{\beta}$ into itself such that $(Y,Z) \in L^p_{\beta}$ is a solution of the BSDE \eqref{S11-eq} if and only if it is a fixed point of $\Phi$.

\noindent More precisely, given $(U,V) \in L^p_{\beta}$, we set $(Y,Z)= \Phi(U,V)$ if  $(Y,Z)$ is the pair satisfying:
\begin{equation}\label{S11-Contr}
Y_t+\int_t^T\int_KZ(s,y)\,(\mu(ds,dx)- \nu(ds,dx))
=\xi+\int_t^T\int_K [f(s,x,U_s,V_s(\cdot)) -V_s(\cdot)]\,\nu(ds,dx).
\end{equation}

\noindent We note that by the Lipschitz condition of $f$ we have
\bee
\begin{array}{l}
\int_t^T \int_K|f(s,x,U_s,V_s)|\,\nu(ds,dx)  = \int_t^T\int_K e^{-\frac{\beta}{p} A_s}e^{\frac{\beta}{p} A_s}|f(s,x,U_s,V_s)|\,\nu(ds,dx) \qquad \qquad \qquad\\
 \,   \leq  \int_t^Te^{-\frac{\beta}{p} A_s}e^{\frac{\beta}{p} A_s}\left[L'|U_s| +L\big(\int_K|V_s|^p\,\phi_s(x)\big)^{1/p}+ \int_K |f(s,x,0,0)|\phi_s(dx) \right]\,dA_s\\
\,   \leq \Big(\int_t^T e^{-\frac{\beta}{p-1}A_s}  dA_s\Big)^{\frac{p-1}{p}} \cdot
\Big[L' \Big(\int_t^Te^{\beta A_s}|U_s|^p dA_s \Big)^{1/p} \\
\quad + L\Big(\int_t^T\int_K e^{\beta A_s}|V_s(x)|^p \,\nu(ds,dx) \Big)^{1/p} + \Big(\int_t^T\int_K e^{\beta A_s}|f(s,x,0,0)|^p \,\nu(ds,dx) \Big)^{1/p}\Big] .
\end{array}
\eee
Since $\frac{\beta}{p-1}\int_t^T e^{-\frac{\beta}{p-1}A_s}  dA_s= e^{-\frac{\beta}{p-1} A_t}-e^{-\frac{\beta}{p-1} A_T}\le e^{-\frac{\beta}{p-1} A_t}$  we
arrive at
\begin{equation}\label{felleunokappabeta}
\begin{array}{l}
\left(\int_t^T\int_K |f(s,x,U_s,V_s)|\,\nu(ds,dx)\, \right)^p \\
 \quad \le C_p\big(\frac{p-1}{\be}\big)^{p-1} e^{-\beta A_t}\int_t^T \int_K e^{\beta A_s}((L')^p|U_s|^p+ L^p|V(s,x)|^p+|f(s,x,0,0)|^p) \,\nu(ds,dx)\\
 \quad \le C_{(p, \beta, L',L)}\int_t^T \int_K e^{\beta A_s}(|U_s|^p+ |V(s,x)|^p+|f(s,x,0,0)|^p) \,\nu(ds,dx).
\end{array}
\end{equation}
Since $(U,V)$ are in $L^p_{\beta}$ and the \eqref{Lp-f0} hold, the last inequality implies in particular that the random variable
 $\int_0^T\int_K|f(s,x,U_s,V_s)|\,\nu(ds,dx)$ is $p$
integrable.

The solution $(Y,Z)$ is defined by considering a c\`{a}dl\`{a}g
version of the martingale $M_t=\E^{\calf_t}\left[\xi + \int_0^T \int_K [f(s,x, U_s,V_s) -V_s]\,\nu(ds,dx)  \right]$. By the martingale representation Theorem
\ref{rappresentazione}, there exists a predictable  process $Z$ with $\E \int_0^T \int_K |Z_s(x)|\, \nu(ds , dx) <\infty$ such that
\begin{equation*} M_t= M_0+ \int_0^t \int_K Z(s,x) \;(\mu(dx \,ds) -\nu(ds,dx)) , \qquad t \in [0,  T].
 \end{equation*}
 Define the process $Y$ by
\begin{equation}\label{Y-M}Y_t = M_t - \int_0^t \int_K [f(s,x,U_s,V_s) -V_s]\,\nu(ds,dx), \qquad t \in [0,  T].
 \end{equation}
Noting that  $Y_T= \xi$, we easily deduce that the equation \eqref{S11-Contr} is satisfied.

It remains to show that $(Y,Z)\in L_{\beta}^p$.
It follows by \eqref{Y-M} that
$$Y_t = \E^{\calf_t}\left[\xi + \int_t^T \int_K [f(s,x,U_s,V_s)-V_s]\,\nu(ds,dx)\right]$$  and so, using \eqref{felleunokappabeta}, we obtain
\begin{eqnarray}\nonumber
e^{\beta A_t} |Y_t|^p  \le
2^{p-1}e^{\beta A_t} |\E^{\calf_t}\xi|^p
+2^{p-1}e^{\beta A_t} \left|\E^{\calf_t}\int_t^T\int_K [f(s,x,U_s,V_s) -V_s]\,\nu(ds,dx) \right|^p \qquad \qquad \qquad\qquad  \qquad \\
 \quad  \quad \le
 2^{p-1}\E^{\calf_t} \left[e^{\beta A_T}|\xi|^p +  C_{(p,\be,L',L,T)}
\int_0^T \int_K  e^{\beta A_s}(|U_s|^p+ |V(s,x)|^p+|f(s,x,0,0)|^p) \,\nu(ds,dx). \right] \qquad \qquad.
 \label{martingausil}
 \end{eqnarray}
Denoting by $m_t$ the right-hand side of \eqref{martingausil}, we see that $m$ is a martingale. In particular, for every stopping time $S$ with values in $[0,T]$,
we have
\begin{equation}\label{uniftempiarresto}
\E\,e^{\beta A_S} |Y_S|^p\le \E\, m_S=\E\,m_T<\infty
\end{equation}
by the optional stopping theorem. Next we define the increasing sequence of stopping times
$$
S_n=\inf\{t\in [0,T]\,:\,
 \int_0^t  \int_Ke^{\beta A_s }( |Y_s|^p+ |Z(s,x)|^p )\,\nu(ds,dx) >n\},
 $$
with the  convention $\inf \emptyset =T$. The It\^{o} formula \eqref{P1} can be applied to $Y$, $Z$
on the interval $[0,S_n]$. Hence, proceeding as in Lemma \ref{LP2}, we deduce
$$
\begin{array}{l}
\E\int_t^T\int_K \left[\frac{1}{2}\left(\be- 1 - \frac{c_{\epsilon}}{1+ \epsilon}\right)|Y_s|^p +\frac{1}{1+ \epsilon}|Z(s,x)|^p\right]\,e^{\be A_s}\,\nu(ds,dx)\\
\qquad \qquad\leq \E\Big(|Y_{S_n}|^pe^{\be A_{S_n}}+
c(\epsilon,\beta,p) \int_t^T\int_K|f(s,x,U_s,V_s)-V_s|^p\,e^{\be A_s}\,\nu(ds,dx)\Big).
\end{array}$$

%dis &\le&\E \,e^{\beta A_{S_n}}|Y_{S_n}|^p +(\frac{\beta-2}{2})\E \int_0^{S_n}\int_K  e^{\beta A_s} |Y_{s}|^p \, \nu(ds,dx)\\
%& + &C(\beta)p^p\E \int_K\int_0^{S_n}  e^{\beta A_s}  |f(s,x,U_s,V_s)|^p\, \nu(ds,dx).
 % \dis &\le& \E \,e^{\beta A_{S_n}}|Y_{S_n}|^p +2\E \int_0^{S_n}\int_K  e^{\beta A_s} |Y_{s}|^p \, \nu(ds,dx)\\
  %& + &p\E \int_0^{S_n}  e^{\beta A_s} |Y_{s}|^{p-1} ( L'|U_s| + L|V_s| + |f(s,x,0,0)| )\, \nu(ds,dx) \\
 % \dis &\le&\E \,e^{\beta A_{S_n}}|Y_{S_n}|^p +(2+\epsilon)\E \int_0^{S_n}\int_K  e^{\beta A_s} |Y_{s}|^p \, \nu(ds,dx)\\
 % &+ &C(\epsilon,p)\E \int_0^{S_n}  e^{\beta A_s} (|U_s|^p +|V_s|^p + |f(s,x,0,0)|^p )\, \nu(ds,dx).

From \eqref{uniftempiarresto} (with $S=S_n$) which we deduce
\begin{equation}\label{stimasonesseenne}
\begin{array}{l}\dis
   \E\int_0^{S_n} \int_K  e^{\beta A_s}(| {Y}_s|^p + |{Z}_s(y)|^p ) \nu(dx,ds)
\\\dis\,\le c_1(\beta,\epsilon,p)
\E e^{\beta A_T}|\xi|^p +  c_2(\beta,\epsilon,p,L',L,T)\int_0^T\int_K e^{\beta A_s}(|U_s|^p+ |V(s,x)|^p+|f(s,x,0,0)|^p) \,\nu(ds,dx).
\end{array}
\end{equation}

Setting $S=\lim_nS_n$ we deduce
$$   \E\int_0^{S} \int_K  e^{\beta A_s}(| {Y}_s|^p + |{Z}_s(y)|^p ) \nu(dx,ds) <\infty,\qquad \P-a.s.
$$
which implies $S=T$, $\P$-a.s., by the definition of $S_n$. Letting $n\to\infty$
in \eqref{stimasonesseenne} we conclude that
 $(Y,Z)\in \L_{\beta}^p$.

Finally we prove that the map $\Phi$ is a contraction.
Let $(U^i,V^i)$, $i=1,2$, be elements of $L^p_{\beta}$ and
let $(Y^i,Z^i)= \Phi (U^i,V^i)$. Denote $\overline{Y}=Y^1-Y^2$,
$\overline{Z}=Z^1-Z^2$, $\overline{U}=U^1-U^2$,
$\overline{V}=V^1-V^2$, $\overline{f}_s=f(s,x,U^1_s,V^1_s)
-f(s,x,U^2_s,V^2_s) -\overline{V}$.  Lemma \ref{LP2} applies to
$\overline{Y},\overline{Z}$. Noting that $\overline{Y}_T=0$, we obtain
$$
\begin{array}{l}\dis
\E e^{\beta A_t}| \overline{Y}_t|^p
+ \E\int_t^T \int_K \left(\frac{1}{1+\epsilon}|\overline{Z}(s,x)|^p +\left(\beta- 1 - \frac{c_{\epsilon}}{1+ \epsilon}\right)  |\overline{Y}_s|^p\right) e^{\beta A_s}\,\nu(ds,dx)  \\
\dis\qquad
\le
  p\E \int_t^T \int_K e^{\beta A_s} |\overline{Y}_{s}|^{p-1}\left[|\overline{f}_s| + |\overline{V}(s,x)|\right]\, \nu(ds,dx) ,
  \qquad t\in [0,T].
\end{array}
$$
From the Lipschitz conditions of $f$ and elementary inequalities it follows that
$$
\begin{array}{l}\dis
\E\int_0^T \int_K  \left[\left(\be- 1 - \frac{c_{\epsilon}}{1+ \epsilon} \right) |\overline{Y}_s|^p
+ \frac{1}{1+ \epsilon}|\overline{Z}(s,x)|^p \right] e^{\beta A_s}\nu(ds,dx)
\\
\quad\le
 p (L+1)\E \int_t^T  e^{\beta A_s} |\overline{Y}_{s}|^{p-1} \,
 \left( \int_K |\overline{V}(s,x)|^p \phi_s(dx)\right)^{1/p}
   dA_s
  + pL^{'}\E \int_t^T \int_K e^{\beta A_s} |\overline{Y}_{s}|^{p-1} \,|\overline{U}_s|\,\nu(ds,dx)
\\
\,
\le
\frac{\alpha}{1+ \epsilon}
\E \int_0^T \int_K e^{\beta A_s}|\overline{V}(s,x)|^p \,\nu(ds,dx)+
(p-1)\left((L+1)^p\cdot\frac{1+ \epsilon}{\alpha}\right)^{\frac{1}{p-1}} \E \int_0^T \int_K e^{\beta A_s}|\overline{Y}_s |^p \, \,\nu(ds,dx)
\\
\quad+ (p-1)\gamma L^{'} \E \int_0^T \int_K e^{\beta A_s}|\overline{Y}_s |^p \,\nu(ds,dx)+
 L^{'}\left(\frac{1}{\gamma}\right)^{p-1}  \E \int_0^T \int_K e^{\beta A_s}|\overline{U}_s |^p \,\nu(ds,dx)
\end{array}
$$
for every $\alpha>0$, $\gamma>0$. This can be written
$$
\begin{array}{l}
\E \int_0^T \int_K \left[\left(\beta  - 1 - \frac{c_{\epsilon}}{1+ \epsilon}- (p-1)\left(\frac{(L+1)^p(1+ \epsilon)}{\alpha}\right)^{\frac{1}{p-1}} - (p-1)\gamma L^{'} \right)\,|\overline{Y}_s|^p \right]\nu(ds,dx)\\
 \qquad + \E \int_0^T \int_K \left[ \frac{1}{1+ \epsilon}|\overline{Z}(s,x)|^p\right] e^{\beta A_s}\,\nu(ds,dx)  \\
\qquad \le \E \int_0^T \int_K
\left[ L^{'}\left(\frac{1}{\gamma}\right)^{p-1}   |\overline{U}_s |^p + \frac{\alpha}{1+ \epsilon}
|\overline{V}(s,x)|^p \right]e^{\beta A_s}\,\nu(ds,dx).
\end{array}
$$
By the assumption on $\beta$ it is possible to find $\alpha\in (0,1)$ such that
$$
\beta >1 + \frac{c_{\epsilon}}{1+ \epsilon}+ (p-1)\left(\frac{(L+1)^p(1+ \epsilon)}{\alpha}\right)^{\frac{1}{p-1}} +{\frac{pL^{'}}{\sqrt[p]\alpha}}.
$$
If $L^{'}=0$ we see that $\Phi$   is an $\alpha$-contraction on
$L^p_{\beta}$ endowed with the equivalent norm
\begin{multline}
(Y,Z)\mapsto
\E \int_0^T \int_K \left[\beta  - 1 - \frac{c_{\epsilon}}{1+ \epsilon}- (p-1)\left(\frac{(L+1)^p (1+ \epsilon)}{\alpha}\right)^{\frac{1}{p-1}}\right]|\overline{Y}_s|^p \nu(ds,dx)\\
+ \E \int_0^T \int_K \left[\frac{1}{1+ \epsilon}|\overline{Z}(s,x)|^p\right] e^{\beta A_s}\,\nu(ds,dx).
\end{multline}
 If $L^{'}>0$ we choose $\gamma=1/\sqrt[p]{\alpha}$
and obtain
$$
\begin{array}{l}
\E \int_0^T \int_K \left[\frac{L^{'}}{\sqrt[p]{\alpha}} |\overline{Y}_s|^p, \nu(ds,dx)
+ \frac{1}{1+ \epsilon}|\overline{Z}(s,x)|^p\right] e^{\beta A_s}\,\nu(ds,dx) \\
\qquad \le
\E \int_0^T \int_K \left[
 L^{'}(\sqrt[p]{\alpha})^{p-1}   |\overline{U}_s  |^p +\frac{1}{1+ \epsilon} \alpha|\overline{V}(s,x)|^p \right] e^{\beta A_s}\,\nu(ds,dx)\\
 \qquad =
 \alpha \,
 \E \int_0^T \int_K \left[ \frac{L^{'}}{\sqrt[p]{\alpha}}|\overline{U}_s|^p
+ \frac{1}{1+ \epsilon}|\overline{V}(s,x)|^p \right] e^{\beta A_s}\,\nu(ds,dx),
\end{array}
$$
so that
 $\Phi$   is an $\alpha$-contraction on
$L^p_{\beta}$ endowed with the equivalent norm
$$(Y,Z)\mapsto \E \int_0^T \int_K \left[\frac{L^{'}}{\sqrt[p]{\alpha}} |{Y}|^p
+\frac{1}{1+ \epsilon}|{Z}|^p\right] e^{\beta A_s}\,\nu(ds,dx).$$ In all cases
there exists a unique fixed point which is the required unique solution to the BSDE \eqref{S11-eq}.
\qed

\begin{remark} \em{Under Hypothesis \ref{hyp:BSDE} we have existence of the solution to the BSDE \eqref{S11}, in the sense of Definition \ref{def:solBSDE}. In contrast, the uniqueness   holds  in the smaller subclasses $L^p_{\beta}$ but it is not guaranteed within the class of all possible solutions  as show the following example.

\noindent Consider a univariate point process. The space $K=\{\De\}$ is
a singleton, and
$N_t~=~1_{\{S\leq t\}},$
where $S$ is a variable with values in $(0,T]\cup\{\infty\}$. The filtration
$(\f_t)$ is still the one generated by $N$, and $G$ denotes the law of $S$,
whereas $g(t)=G((t,\infty]$.
We suppose that $G$ has
no atom, but is supported by $[0,T]$. We have $A_t=a(t\wedge S)$, where
$a(t)=-\log g(t)$ is increasing, finite for $t<v$ and infinite if
$t\geq v$, where $v=\inf(t:\,g(t)=0)\leq T$ is the right end point of the
support of the measure $G$.

\noindent Consider the following
 equation with $\xi=0$ and $f(t,x,y,z)=z$
\bee\label{B21}
Y_t+\int_{(t,T]}Z_s\,(dN_s-dA_s)=0,
\eee

\noindent  $Y_t=0$ is solution, but $Y_t=we^{A_t}\,1_{\{t<S\}}$ for any $w\in\R$ is also
a solution (see \cite[Proposition 11]{CFJ})
.}
\end{remark}

 \section{An approximation scheme}\label{S4}
In this section we show how it is possible to
reduce the problem of solving Equation \eqref{S11} to solving   a finite system of ordinary differential equation: the solution to the \eqref{S11} can be obtained as limit of a finite sequence of deterministic differential equation.

 \subsection{An approximation of the BSDE \eqref{S11}.} We will approximate the BSDE \eqref{S11} by a other BSDE driven by random measures with a finite number of jumps.

 \noindent For each (finite) integer $m\ge 1$ let us consider
 the BSDE
 \bee\label{C3}
\begin{array}{c}
Y^m_t+\int_t^T\int_EZ^m(s,x)\,\mu^m(ds,dx)
=\xi^m+\int_t^T\int_Ef(s,x,Y_{s}^m,Z^m_s(\cdot))\,1_{s\le S_{m}}\,
\nu^m(ds,dx) \\
\mu^m(ds,dx)=\mu(ds,dx)\,1_{\{s\leq S_m\}},\quad
\nu^m(ds,dx)=\nu(ds,dx)\,1_{\{s\leq S_m\}},
\quad\xi^m=\xi\,1_{\{T<S_m\}}.
\end{array}
\eee
Then $\nu^m$ is the compensator of $\mu^m$, relative to $(\f_t)$ and
also to the smaller filtration $(\f^{(m)}_t=\f_{t\wedge S_m})$ generated by
$\mu^m$, whereas $\xi^m$ is $\f^{(m)}_T$-measurable.

\noindent We note that the
generator $f(s,x,y,z)\,1_{s\le S_{m}}$ and the terminal condition
$\xi^m$ satisfy Hypothesis \ref{hyp:BSDE}
with respect to    $(\calf_{t\wedge S_m})$.
By Theorem \ref{Thm:existence!bsde}  there exists a solution
$(Y^m,Z^m)$ in $L^p_{\beta}$ for $\beta>1 + \frac{c_{\epsilon}}{1+ \epsilon}+ pL'+(p-1)((L+1)^p(1+ \epsilon))^{\frac{1}{p-1}}$.
We may assume
$Y^m_s= \xi^m =0,\;Z^m_s=0$ for $ S_m< s\le T$, so that the BSDE \eqref{C3} becomes
\begin{equation}\label{BSDEtruncated}
    Y^m_t+\int_t^T\int_K Z^m(s,x)\,\mu(ds,dx)=\xi^m
    +\int_t^T \int_Kf(s,x,Y^m_s,Z^m(s,x))\,1_{s\le S_{m}}\,\nu(ds,dx),
\qquad t\in [0,T],
\end{equation}
and we have
$
    \E\left[
    \int_0^T\int_K(|Y_s^m |^p+|Z(s,x) ^m|^p)\,e^{\beta A_s}\nu(ds,dx) \right]<\infty.$
    Clearly, $Y^m$ is  adapted and $Z^m$ is predictable
 also   with respect to the natural filtration $(\calf_t)$ of $\mu$.
 \begin{proposition}\label{approxy,z}
Let $(Y^m,Z^m)$  and and $(Y,Z)$ be the  solutions in $L^p_{\beta}$ to the BSDEs \eqref{BSDEtruncated} and  \eqref{S11} respectively. Then, if  $\beta>1 + \frac{c_{\epsilon}}{1+ \epsilon}+ pL'+(p-1)((L+1)^p(1+ \epsilon))^{\frac{1}{p-1}}$, there  exists a constant $C$ independent of $m$ such that

$$\begin{array}{l}\dis
  \E
    \int_0^T\int_K(|Y_s-  Y^m_s |^p+|Z(s,x) -Z^{m}(s,x)|^p)\,e^{\beta A_s}\nu(ds,dx)
        \\\dis
    \le
   C \E\left[ |  \xi|^pe^{\beta A_T}1_{S_m<T}+
    \int_{S_m\wedge T}^{  T} \int_K |f(s,x,0,0)|^pe^{\beta A_s}\nu(ds,dx) \right]\\

\end{array}
    $$
\end{proposition}

{\bf Proof.}

We define
$\bar Y=Y-Y^{m}$,
$\bar Z=Z-Z^m$,
and note that $(\bar Y,\bar Z)$ is the solution to
$$
    \bar Y _t+\int_t^T \int_K \bar Z(s,x)\,\mu(ds,dx)=\bar \xi +\int_t^T\int_K
    \bar f(s,x,\bar Y_s,\bar Z(s,x))\,\nu(ds,dx)
\qquad t\in [0,T],
$$
where we have set $\bar f(s,x,y,z)=
f(s,x,Y_s^m+y,Z_s^m+z)- f(s,x,Y_s^m,Z_s^m)1_{s\le S_{m}}$
 and $\bar \xi = \xi-\xi^m$.
We deduce from the a priori estimate
\eqref{P4}  that
$$
\begin{array}{l}\dis
\,    \E
    \int_0^T\int_K(|\bar Y_s |^p+|\bar Z(s,x) |^p)\,e^{\beta A_s} \nu(ds,dx)
        \\\dis
    \le
    C\E\left[ |\bar \xi|^p e^{\beta A_T}+
    \int_0^T|\bar f(s,x,0,0)|^p e^{\beta A_s}\nu(ds,dx) \right]
        \\\dis
    =
    \E\left[ |  \xi|^pe^{\beta A_T}1_{S_m<T}+
    \int_{S_m\wedge T}^{ T} \int_K |f(s,x,0,0)|^p e^{\beta A_s}\nu(ds,dx) \right]
\end{array}
    $$
where the last equality holds because  $Y^m_s=0,\;Z^m_s=0$ for $ S_m< s\le T$.
It follows that
\bee\label{error}
\begin{array}{l}\dis
  \E
    \int_0^T\int_K(|Y_s-  Y^m_s |^p+|Z(s,x) -Z^{m}(s,x)|^p)\,e^{\beta A_s}\nu(ds,dx)
        \\\dis
    \le
   C \E\left[ |  \xi|^pe^{\beta A_T}1_{S_m<T}+
    \int_{S_m\wedge T}^{  T} \int_K |f(s,x,0,0)|^pe^{\beta A_s}\nu(ds,dx) \right].
\end{array}
\eee
\qed

\begin{remark} \em Since the right-hand side in \eqref{error} tends to $0$ as $m\to\infty$, Proposition \ref{approxy,z} ensures that it is possible approximate the solution $(Y,Z)$ to the \eqref{S11}, with the solution to the BSDE \eqref{BSDEtruncated} driven by random measures with a finite number of jumps. Moreover it furnish the error estimate for this approximation. In the case when $\xi$ and $f(s,x,0,0)$ are uniformly bounded, the approximation error can easily be expressed in terms of $P(S_m <T)$. Under  $L^p_{\beta}$ integrability conditions in Hypothesis \ref{hyp:BSDE}, a similar result can be obtained using  H\"{o}lder inequality.
\end{remark}

\subsection{Reduction to ordinary differential equations.} \label{subsect: ode}
In this section, we recall how, using a result in \cite{CFJ}, it is possible
to reduce the problem of solving equation \eqref{BSDEtruncated} to solving a sequence of ordinary
differential equations.

\noindent Beside $(\bf{A_1})$ we need here also the following

$ $

\noindent \textbf{Assumption $\bf{A_2}$ :} \rm $\P(S_{n+1}>T\mid\f_{S_n})>0$ for all $n\ge 0$.

$ $

This condition is useful to characterize the
$\f_{S_n}$-conditional law of $(S_{n+1},X_{n+1})$ and the compensator $\nu$
of $\mu$ (see \cite[Section 4]{CFJ}).

The process $(S_n,X_n)$ takes its values in the set $\s=([0,T]\times K)
\cup\{(\infty,\De)\}$.
For any integer $n\geq0$ we let $H_n$ be the subset of $\s^{n+1}$
consisting in all $D=((t_0,x_0),\cdots,(t_n,x_n))$ satisfying
$$t_0=0,~x_0=\De,\quad t_{j+1}\geq t_j,\quad
t_j\leq T~\Rightarrow~t_{j+1}>t_j,\quad t_j>T~\Rightarrow~
(t_j,x_j)=(\infty,\De).$$
We set $D^{\max}=t_n$ and endow $H_n$ with its Borel $\si$-field $\h_n$.
We set $S_0=0$ and $X_0=\De$, so
\bee\label{B2}
D_n=((S_0,X_0),\cdots,(S_n,X_n))
\eee
is a random element with values in $H_n$, whose law is denoted as $\La_n$
(a probability measure on $(H_n,\h_n)$).

The process $Y^m$ solution to \eqref{BSDEtruncated} is an adapted c\`adl\`ag process, which
further is continuous outside the times $S_n$. Hence for
each $ 0 \leq n\le m$ there is a Borel function $y^n=y^n_D(t)$ on $H_n\times[0,T]$ such
that
\bee\label{B4}
\begin{array}{l}
D^{\max}=\infty~~\Rightarrow~~y^n_D(t)=0\\
t\mapsto y^n_D(t)~~\text{is continuous on $[0,T]$ and constant on
$[0,T\wedge D^{\max}]$}\\
S_n(\om)\leq t< S_{n+1}(\om),~~t\leq T~~\Rightarrow~Y_t(\om)=y_{D_n(\om)}^n(t),
\end{array}
\eee
and we express this as $Y\equiv(y^n)_{n=0}^m$. Also the component $Z^m$ of the solution to \eqref{BSDEtruncated} can be express  as $Z^m\equiv(z^n)_{n=0}^m$
 where
$z^n=z^n_D(t,x)$ is a Borel function on $H_n\times[0,T]\times E$ such that
\bee\label{B5}
\begin{array}{l}
D^{\max}=\infty~~\Rightarrow~~z^n_D(t,x)=0\\
S_n(\om)<t\leq S_{n+1}(\om)\wedge T~~\Rightarrow~~
Z(\om,t,x)=z_{D_n(\om)}^n(t,x).
\end{array}
\eee
The generator $f\,1_{s\le S_{m}}$ has a nice predictability
property only after plugging in a predictable function $Z$. This implies
that, for any $0 \leq n \leq m$, and  $z^n=z^n_D(t,x)$, one has
a Borel function $f\{z^n\}^{n}=f\{z^n\}^{n}_D(t,x,y,w)$ on
$H_n\times[0,T]\times E\times\R\times\R$, such that (with $t\leq T$ below)
\bee\label{B8}
\begin{array}{l}
D^{\max}=\infty~~\Rightarrow~~f\{z^n\}^n_{D}(t,x,y)=0\\
S_n(\om)<t\leq S_{n+1}(\om),~\ze(x)=w+z^n_{D_n(\om)}(t,x)
~\Rightarrow~
f(\om,t,x,y,\ze)=f\{z^n\}^n_{D_n(\om)}(t,x,y,w).
\end{array}
\eee

The variable $\xi^m$ is $\f_T$-measurable, hence
for each $0\leq n \le m$ there is an $\h_n$-measurable map $D\mapsto u^n_D$ on $H_n$
with
\bee\label{B3}
\begin{array}{l}
D^{\max}=\infty~~\Rightarrow~~u^n_D=0\\
S_n(\om)\leq T<S_{n+1}(\om)~~\Rightarrow~~\xi(\om)=u_{D_n(\om)}^n.
\end{array}
\eee

\noindent By \cite[Lemma 7]{CFJ}
 $Y^m\equiv(y^n)_{n=0}^m$ is a solution if and
only if for $P$-almost all $\om$ we have:
\bee\label{B14}
t\in[0,T]~~\Rightarrow~~y^m_{D_m(\om)}(t)=u^m_{D_m(\om)}=\xi(\om).
\eee
and  for all
$n=0,\cdots,m-1$
\bee\label{B13}
y^n_{D_n(\om)}(t)=u^n_{D_n(\om)}+\int_t^T\int_E f\{\widehat{y}^{n+1}\}^n_{D_n(\om)}(s,x,y^n_{D_n(\om)}(s),
-y^n_{D_n(\om)}(s))\,\nu^n_{D_n(\om)}(ds,dx),
\, t\in[0,T].
\eee
where we set
\bee\label{B10}
\widehat{y}^{n+1}=(\widehat{y}^{n+1}_D(t,x):\,(D,t,x)\in H_n\times[0,T]\times E):\quad
\widehat{y}^{n+1}_D(t,x)=y_{D\cup\{(t,x)\}}^{n+1}(t)\,1_{\{t>D^{max}\}}.
\eee

As stressed  before, we need to assume $(\bf{A_1})$, as in previous sections, and also $(\bf{A_2})$. The following Lemma provides a condition which implies Assumption $(\bf{A_2})$.
\begin{lem}\label{lemmacondprobabpos}
Assume $(\bf{A_1})$. Then, for  $n\ge 0$
 the inequality $\E [e^{A_{T\wedge S_{n+1}}}]<\infty$ implies
\begin{equation}\label{condprobabpos}
    \P(S_{n+1}>T\,|\,\calf_{S_n})>0 \quad a.s.
\end{equation}

In particular, if
\begin{equation}\label{eAT}\E[e^{\beta A_T}] < \infty,
\end{equation} then
\eqref{condprobabpos} holds true for every $n\ge 0$
.
\end{lem}

{\bf Proof.} Let
$G'^n_{D_n}(dt)$ be the conditional law of $S_{n+1}$ given $\calf_{S_n}$.
Let us introduce the corresponding cumulative distribution function
$F_D(t)=G'^n_D((0,t])$.
Since we assume that the
dual predictable projection $A$ of $\mu$ is continuous we can take a version
of $F_D$ which is continuous in $t$ and we have, $\P$-a.s.,
\begin{equation}\label{logdistrfunct}
    \begin{array}{lll}
A_t&=&\dis
A_{S_n}+\int_{S_n}^t\frac{1}{1-F_{D_n}(s)}\,F_{D_n}(ds) =
\\
&=&\dis A_{S_n}-\log (1-F_{D_n}(t)),
\qquad S_n< t\le S_{n+1}.
\end{array}
\end{equation}
Since $F_D$ is continuous in $t$, the  conditional law of
$F_{D_n}(S_{n+1})$ given $\calf_{S_n}$ is the uniform distribution
on $(0,1)$, so that in particular
$\E [(1-F_{D_n}(S_{n+1}))^{-1}|\calf_{S_n}]=\infty$ a.s.
Now suppose that \eqref{condprobabpos} is violated for some $n$. Then
there exists $Q\in \calf_{S_n}$ with $\P(Q)>0$ such that
$    \P(S_{n+1}\le T\,|\,\calf_{S_n})=1$ on $Q$.  Then
$$
\P(Q)= \E[1_Q\, \P(S_{n+1}\le T\,|\,\calf_{S_n})]
= \P(Q\cap \{S_{n+1}\le T\}),
$$
which shows that $S_{n+1}\le T$ a.s. on $Q$. It follows from
\eqref{logdistrfunct}
that
$$
\begin{array}{l}
\dis
\E [e^{  A_{T\wedge S_{n+1}}}]\ge\E [1_Q\,e^{  A_{T\wedge S_{n+1}}}]\ge
\E[1_Q\,  e^{  A_{S_{n+1}}}] \ge
\\\dis
\E[1_Q\,(1-F_{D_n}(S_{n+1}))^{-1}]=
\E[1_Q\,\E [(1-F_{D_n}(S_{n+1}))^{-1}|\calf_{S_n}]]=\infty,
\end{array}
$$
contradicting the assumption. The first part of the lemma is therefore proved.

Next assume the \eqref{eAT}. Then, the conclusion follows from the statement proved above noting that
$$
\begin{array}{l}
\dis
\infty>
\E[e^{\beta A_T}]\ge
\E[ e^{\beta A_{T\wedge S_n}}].
\end{array}
$$
\qed

\section{Optimal control}\label{S5}

In this section we use the previous results on the BSDE to solve an optimal control problem.
We assume that  a marked point process $\mu$ of \eqref{S1} is given on $(\Omega, \mathcal{F}, \mathbb{P})$, generating the filtration $(\mathcal{F}_t)$ and satisfying $\bf{A_1}$.
In particular we suppose that $T_n\to\infty$ $\P$-a.s.

The data specifying the optimal control problem are an
  action (or decision) space
$U$, and a function $r$ specifying the effect of the control
process, a running cost function $l$, a terminal cost function $g$. We assume that these data  satisfy the following conditions.

\begin{hypothesis}\label{hyp:controllo}
\begin{enumerate}
\item  $(U,\calu)$ is a measurable space.
\item  The functions $r,l:\Omega\times [0,T]\times K\times U\to \R$
are predictable and
 there exist
constants $C_r> 1$, $C_l>0$ such that
\begin{equation}\label{ellelimitato}
0\le r (\omega,t,x,u)\le C_r,\quad |l(\omega,t,x,u)|\le C_l,\qquad  \omega \in \Om, t\in [0,T],
x\in K, u\in U.
\end{equation}
\item
The function $g:\Omega\times K\to \R$ is $\calf_T\otimes
\calk$-measurable.
\end{enumerate}
\end{hypothesis}

We define as an admissible  control process
any predictable process $(u_t)_{t\in[0,T]}$ with values in $U$. The set of admissible
control processes is denoted $\cala$.

To every control $u(\cdot)\in\cala$ we associate a probability measure $\P_u$ on
$(\Omega,\calf)$ by a change of measure of Girsanov type, as we now describe. We define
$$
L_t=
\exp\left(\int_0^t\int_K (1-r (s,x,u_s))\;\nu(ds,dx)\right)
\prod_{n\ge1\,:\,S_n\le t}r(S_n,X_n,u_{S_n}),
\qquad t\in [0,T],
$$
with the convention that the last product equals $1$ if there are
no indices $n\ge 1$ satisfying $S_n\le t$.  $L$ is a nonnegative supermartingale, (see \cite{J74}
Proposition 4.3). Moreover the following result holds

\begin{lem}\label{Con-Fuh}(\cite[Lemma 4.2]{CoFu-mpp})
Let $\gamma>1$ and
\begin{equation}\label{valoredibeta}
    \beta= \gamma+1+\frac{C_r^{q^2}}{q-1}.
\end{equation}
If $\E \exp (\beta A_T)<\infty$ then  we have $\sup_{t\in [0,T]}\E L_t^{\gamma}<\infty$ and
$\E L_T=1$.
\end{lem}
Under the assumption of the lemma, the process $L$ is a martingale.
By Girsanov's Theorem for point
processes, the predictable compensator of the measure $\mu$ under $\P_u$ is
$$\nu^u(dt,dx)=r(t,x,u_t)\,\nu(dt,\,dx)=r(t,x,u_t)\,\phi_t(dx)\,dA_t.$$
We finally define the cost associated to every $u(\cdot)\in\aaa$ as
$$J(u(\cdot))=\E_u\Big(\int_0^Tl(t,X_t,u_t)\,dA_t + g(X_T) \Big),$$
where $\E_u$ denotes the expectation under $\P_u$.
Later we will
assume that
\begin{equation}
\label{condsucostofinale}
 \E[|g(X_T)|^pe^{\beta A_T} ]<\infty
 \end{equation}
  for some $\beta>0$
that will be fixed in such a way that the cost is finite for every
admissible control. The control problem consists in minimizing $
J(u(\cdot))$ over $\cala$.
A basic role is played by the BSDE
\begin{equation}\label{bsdecontrollo}
Y_t+\int_t^T\int_KZ(s,x)\,\mu(ds,dx) =
\xi  +\int_t^T f(s,X_s,Z_s(\cdot))\,dA_s,
\end{equation}
with terminal condition $g(X_T)$ being the terminal cost above, and with the
generator $f$ being the Hamiltonian function defined below. This is Equation
\eqref{S11}, with $f$ only depending on $(\om,t,\ze)$, and indeed it comes
from an equation of type II via the transformation \eqref{S14}.

The Hamiltonian function $f$ is defined on $\Om\times[0,T]\times\ba(E)$ as
\bee\label{CC6}
f(\om,t,x,z(\cdot))=\left\{\begin{array}{ll}
\inf_{u\in U}\Big(l(\om,t,x,u)+ \int_K z(x)\,r(\om,t,x,u)\,\phi_{\om,t}(dx)\Big)
&\text{if}~\int_K|\ze(x)|\,\phi_{\om,t}(dx)<\infty\\
0&\text{otherwise.}\end{array}\right.
\end{equation}

We will assume that the infimum is in fact achieved, possibly at
many points. Moreover we need to verify that the generator of the
BSDE satisfies the conditions required in Hypothesis \ref{hyp:BSDE}, in particular
the measurability property which does not
follow from its definition.
An appropriate assumption is the following one, since we
will see below in Proposition \ref{PCC1} that it can
be verified under quite general conditions.
\vsc

\begin{hypothesis}\label{exist-min-ham}
For every predictable function $Z$ on $\Omega\times [0,T]\times E$
there exists a $U$-valued predictable process
(i.e., an admissible control) $\underline{u}^Z$ such that,
$dA_t(\omega)\P(d\omega)$-almost surely,
\begin{equation}\label{minselector}
    f(\omega,t,X_{t-}(\omega),Z_{\omega,t}(\cdot))=
l(t,X_{t-}(\omega),\underline{u}^Z(\omega,t )) + \int_K Z_{\omega,t}(x) \,
\Big(r
(\omega,t,x,\underline{u}^Z(\omega,t))\Big)\,\phi_{\omega,t}(dx)
\end{equation}
\end{hypothesis}

We can now verify that all the assumptions of Hypothesis
\ref{hyp:BSDE} hold true for the generator of the BSDE
\eqref{bsdecontrollo}. The \eqref{minselector} shows that the process $(\omega,t)\mapsto
f(\omega,t,X_{t-}(\omega),Z_t(\omega,\cdot))$ is progressive;
since $A$ is assumed to have continuous trajectories and $X$ has piecewise
constant paths, the progressive set
$\{(\omega,t)\,:\, X_{t-}(\omega)\neq X_{t}(\omega)\}$ has
measure zero with respect to $dA_t(\omega)\P(d\omega)$;  it follows that the
process
$$(\omega,t)\mapsto f(\omega,t,X_t(\omega),Z_t(\omega,\cdot))
$$
is progressive, after modification on a set of measure zero,
as required in \eqref{fmisurabile}.
Using
the boundedness assumptions \eqref{ellelimitato},
 it is easy to check  that \eqref{hyp:f-lip} is verified with
$L'=0$ and
$L=C_r.$
Using \eqref{ellelimitato} again we also have
\begin{equation}\label{gensommcontrollo}
  \E \int_0^T e^{\beta A_t}|f(t,X_t,0)|^p dA_t
  = \E \int_0^T e^{\beta A_t}|\inf_{u\in U}  l(\om,t, X_t,u)|^p dA_t
  \le C_l^p\,\beta^{-1}\,
\E \,e^{\beta A_T},
\end{equation}
 so  that \eqref{Lp-f0} holds as well provided the
 right-hand side of \eqref{gensommcontrollo}
 is finite.
 Assuming finally that \eqref{condsucostofinale} holds, by
Theorem \ref{Thm:existence!bsde} the BSDE has a unique solution
$(Y,Z)\in L^p_{\beta}$ if $\beta>1 + \frac{c_{\epsilon}}{1+ \epsilon}+(p-1)(2L^p(1+ \epsilon))^{\frac{1}{p-1}}$.

The corresponding admissible control $\underline{u}^Z$, whose
existence is required in Assumption (C), will be denoted as $u^*$.

\begin{theorem}\label{TCC1} Assume that Hypotheses \ref{hyp:controllo} and \ref{exist-min-ham} are satisfied and that .
\begin{equation}\label{ATintegrabile}
\E \exp \left[ \left(\frac{2p-1}{p-1}+(p-1)C_r^{\left(\frac{p}{p-1} \right)^2}\right) A_T\right]<\infty.
\end{equation}
Suppose also that there exists $\beta$ such that
\begin{equation}\label{compatibbeta}
  \beta>1 + \frac{c_{\epsilon}}{1+ \epsilon}+(p-1)(({C_r}+1)^p(1+ \epsilon))^{\frac{1}{p-1}}, \quad \E \exp \left(\beta A_T\right)<\infty,
\quad \E[|g(X_T)|^pe^{\beta A_T} ]<\infty.
\end{equation}
Let $(Y,Z)\in L^p_{\beta}$ denote
 the solution to the BSDE \eqref{bsdecontrollo} and
$u^*=\underline{u}^Z$ the corresponding admissible control.
Then
$u^*(\cdot)$ is   optimal and
$Y_0= J(u^*(\cdot))= \inf_{u(\cdot)\in\aaa }J(u(\cdot))$ is the
optimal cost.
\end{theorem}

 \noindent{\bf Proof.} Fix $u(\cdot)\in\cala$.
Assumption \eqref{ATintegrabile} allows to apply Lemma
\ref{Con-Fuh} with $\gamma= \frac{p}{p-1}$ and yields $\E L_T^{\frac{p}{p-1}}<\infty$.
 It follows that $g(X_T)$ is integrable under $\P_u$. Indeed by \eqref{condsucostofinale}
 $$
\E_u |g(X_T)|= \E |L_Tg(X_T)|\le (\E L_T^{\frac{p}{p-1}})^{\frac{p-1}{p}}(\E g(X_T)^p)^{1/p}<\infty.
$$

 We next  show that under $\P_u$ we have
 $\E_u \int_{0}^T \int_K |Z(t,x)|\,\nu^u(dt\,dx)
<\infty$.
First note that, by H\"older's inequality,
$$
\begin{array}{lll}\dis
\int_{0}^T \int_K  |Z(t,x)|\,\nu(dt,dx) &=&\dis
  \int_{0}^T \int_K e^{-\frac{\beta}{p}A_t}e^{\frac{\beta}{p}A_t} |Z(t,x)|\,\nu(dt,dx)
\\
&\le&\dis
\left( \int_{0}^T   e^{-\frac{1}{p-1}\beta A_t} dA_t\right)^{\frac{p-1}{p}}
\left(
\int_{0}^T \int_K e^{\beta A_t} |Z(t,x)|^p\,\nu(dt,dx)
 \right)^{1/p}
\\
&\le&\dis
\left( \frac{p-1}{\beta} \right)^\frac{p-1}{p}
\left(
\int_{0}^T \int_K e^{\beta A_t} |Z(t,x)|^p\,\nu(dt,dx)
 \right)^{1/p}.
\end{array}
$$
Therefore, using \eqref{ellelimitato},
$$
\begin{array}{lll}\dis
\E_u \int_{0}^T \int_K |Z(t,x)|\,\nu^u(dt,dx)&=&\dis
\E_u \int_{0}^T \int_K   |Z(t,x)|\, r(t,x,u_t)\,\nu(dt,dx)
\\&
=&\dis
\E\left[L_T \int_{0}^T \int_K |Z(t,x)|\, r(t,x,u_t)\,\nu(dt,dx)\right]
\\&
\le&\dis
(\E L_T^{\frac{p}{p-1}})^{\frac{p-1}{p}}C_r\left( \frac{p-1}{\beta} \right)^\frac{p-1}{p}
\left(\E\int_{0}^T \int_K  e^{\beta A_t}|Z(t,x)|^p\,\nu(dt,dx)\right)^{1/p}
\end{array}
$$
and the right-hand side of the last inequality is finite,
since $(Y,Z)\in L^p_{\beta}$.

By similar arguments we also check that
$$\begin{array}{lll}
\E_u\int_{0}^T |f(t,Z_t(\cdot))|\,dA_t
&=&\E L_T \,\int_{0}^T  |f(t,X_t,Z_t(\cdot))|\,dA_t \\
&\le &\E L_T \, C_r\Big(\int_{0}^T
\left[\left(\int_K |Z(t,x)|^p\,\phi_{\omega,t}(dx)\right)^{1/p}+|f(t,X_t,0)|\right]\,dA_t\Big)<\infty.
\end{array}$$
Setting $t=0$ and taking the $\P_u$-expectation in the BSDE (\ref{bsdecontrollo})
we therefore obtain
$$Y_0+\E_u\Big(\int_{0}^T \int_K Z(t,x)\, r(t,x,u_t) \,\nu(dt,dx)\Big)=
\E_u(g(X_T))+\E_u\Big(\int_0^T f(t,Z_t(\cdot))\,dA_t\Big).$$
Adding $\E_u\big(\int_{0}^T l(t,X_t,u_t)\,dA_t\big)$
to both sides, we finally obtain the equality
$$\begin{array}{l}
Y_0+\E_u\Big(\int_{0}^T \big(l(t,X_t,u_t)+\int_K Z(t,x)\,r(t,x,u_t)\,\phi_t(dx)
\big)\,dA_t\Big)\\\qquad
= J(u(\cdot))+\E_u\Big(\int_{0}^T f(t,Z_t(\cdot))\,dA_t\Big)\\\qquad
=J(u(\cdot))+\E_u\Big(\int_{0}^T\inf_{u\in U}\big(l(t,X_t,u)+\int_K
Z(t,x)\,r(t,x,u_t),\phi_t(dx)\big)\,dA_t\Big).
\end{array}$$
This implies
immediately the inequality $Y_0 \le  J(u(\cdot))$
for every admissible control, with an equality if $u(\cdot)=u^*(\cdot)$.\qed
\vsc

Assumption (C) can be verified in specific
situations when it is possible to compute explicitly the function
$\underline{u}^Z$. General conditions for its validity can also be
formulated using appropriate measurable selection theorems, as in  the
following proposition. Its proof can be found in \cite{CFJ}, Proposition 17.

\begin{proposition} \label{PCC1} Suppose that $U$ is a compact metric space with its
Borel $\sigma$-field $\ua$ and that the functions
$r(\omega,t,x,\cdot),l(\omega,t,\cdot)$  are continuous on $U$ for every
$(\omega,t,x)$. Then if further $r$ satisfies \eqref{ellelimitato}, Assumption (C) holds.
\end{proposition}

\textbf{Acknowledgments.} We wish to thank Prof. Jean Jacod for discussions on connections between BSDEs
and point processes.

\end{document}